\documentstyle[amscd, amstex]{amsart}

\def\ZZ         {{\Bbb Z}}
\def\RR         {{\Bbb R}}
\def\CC         {{\Bbb C}}
\def\QQ         {{\Bbb Q}}
\def\PP         {{\Bbb P}}

\newtheorem{thm}{Theorem}[section]
\newtheorem{lem}[thm]{Lemma}
\newtheorem{cor}[thm]{Corollary}
\newtheorem{pr}[thm]{Proposition}

\theoremstyle{definition}
\newtheorem{rem}[thm]{Remark}
\newtheorem{ex}[thm]{Example}
\newtheorem{defn}[thm]{Definition}
\newcommand{\st}{{\rm Star}}
\newcommand{\psx}{{{\bf P}_{\Sigma_X}}}

\newcommand{\im}{{\rm im}}
\newcommand{\key}{\bibitem}
\newcommand{\psd}{{{\bf P}_{\Sigma_D}}}
\newcommand{\ps}{{{\bf P}_{\Sigma}}}

\newcommand{\gr}{{\rm Gr}}

\newcommand\hidot{{\raise1pt\hbox{$\scriptscriptstyle\bullet$}}}
\newcommand\lodot{{\raise.3pt\hbox{$\scriptscriptstyle\bullet$}}}
\newcommand{\dd}{{\rm d}}

\begin{document}

\title[Cohomology  of rational forms  on toric varieties]
{Cohomology of rational forms and \\ a vanishing theorem on toric
varieties}
\author{Anvar R. Mavlyutov}
\address {Department of Mathematics, Oklahoma State University, Stillwater, OK 74078, USA.}
 \email{mavlyutov@@math.okstate.edu}

\begin{abstract}
We explicitly describe cohomology of the sheaf of differential
forms with poles along a semiample divisor on a complete
simplicial toric variety. As an application, we obtain a new
vanishing theorem which is an analogue of  the
Bott-Steenbrink-Danilov vanishing theorem.
\end{abstract}

\keywords{Toric varieties, semiample divisors, vanishing theorem.}
\subjclass{Primary: 14M25}

\maketitle

\setcounter{section}{-1}

\section{Introduction.}

 Phillip Griffiths in
\cite{g} calculated the cohomology of smooth hypersurfaces in a
projective space. His method used the Gysin exact sequence, and
the problem was reduced to finding the cohomology of the
complement of the hypersurface in its ambient space. The latter
cohomology was easily found due to the  vanishing theorem of
R.~Bott in \cite{bo}: $$H^k({\Bbb P}^m,\Omega^l_{{\Bbb
P}^m}(X))=0\quad\text{ for an ample divisor }X\text{  and }k>0.$$
This theorem was extended  to  an ample divisor
 on a   complete toric variety (see \cite{bflm}, \cite{d} and
 \cite{bc}).

 As in the case of projective hypersurfaces, in order to compute the cohomology of $X$ quasismooth hypersurfaces in
complete simplicial toric varieties $\ps$, one needs to know the
cohomology of the twisted sheaves
$$H^{k}(\ps,\Omega_{\ps}^{l}(X)).$$ An important case to consider
is when the divisor $X$ is semiample. For toric varieties this
means that the corresponding line bundle is generated by global
sections. The Bott  vanishing theorem does not hold for semiample
divisors. However, in this paper we explicitly  calculated
  the cohomology
  of the twisted sheaves $\Omega_{\ps}^{l}(X)$ and, in particular, we
discovered a new vanishing theorem for semiample divisors:
\begin{equation}\label{e:van}
H^{k}(\ps,\Omega_{\ps}^{l}(X))=0\text{ if }k>l\text{ or }l>k+i,
\end{equation}
where $i$ is the Kodaira-Iitaka dimension of the divisor $X$.
 When
the divisor $X$ is trivial this reduces to the well-known
Danilov's vanishing result for the cohomology of a complete
simplicial toric variety: $$H^k(\ps,\Omega_\ps^l)=0\quad\text{ for
}k\ne l.$$ Moreover, for $k=l$, we recover Danilov-Jurkevich's
description of the cohomology of a complete simplicial toric
variety.

The plan of the paper is as follows. In section~\ref{s:one}, we
compute cohomology of Ishida's complexes that are necessary in
studying cohomology of sheaves $\Omega_{\ps}^{l}(X)$. Then, in
section~\ref{s:two}, we calculated
$H^{k}(\ps,\Omega_{\ps}^{l}(X))$ for a semiample divisor $X$ on a
complete simplicial toric variety $\ps$. As a consequence of this
result, we obtain the vanishing theorem (\ref{e:van}). We also
find a dimension formula  for this cohomology and explicitly
represent the generators of cohomology
$H^{k}(\ps,\Omega_{\ps}^{l}(X))$ in terms of the \v Cech cocycles.

{\it Acknowledgment.} We would like to thank Evgeny Materov for
his interest and for  pointing out a useful reference \cite{is} of
M. Ishida that helped in calculations. We are also grateful to
David Cox for correcting some of the references.

\section{Cohomology of Ishida's complexes.}\label{s:one}

Ishida's complexes have appeared as a result of studying sheaves
of differential forms on toric varieties. In this section, we will
compute the cohomology of Ishida's complexes of modules.

First, we fix some standard notation: $M$ is a lattice of rank
$d$; $N=\text{Hom}(M,{\Bbb Z})$ is the dual lattice; $\Sigma$ is a
finite rational (usually simplicial) fan in the $\Bbb R$-scalar
extension $N_{\Bbb R}$; $\Sigma(k)$ is the set of all
$k$-dimensional cones in $\Sigma$;   $e_1,\dots,e_n$ are the
minimal integral (primitive) generators of  the 1-dimensional
cones $\rho_1,\dots,\rho_n\in\Sigma(1)$.

\begin{defn}\label{d:ish}\cite{o} Let $\Sigma$ be a fan in $N_\RR$ and
$l=0,\dots,d$. Then Ishida's $l$-th complex  of $\QQ$-modules is
denoted $C^*(\Sigma,l)$, where
$$C^j(\Sigma,l)=\bigoplus_{\begin{Sb}\gamma\in\Sigma(j)\end{Sb}}
{\bigwedge}^{l-j}\gamma^\perp $$ (for simplicity, $\gamma^\perp:=
M_\QQ\cap\gamma^\perp$) and the coboundary homomorphism $\delta$
is defined as the direct sum of
$\delta_{\gamma,\tau}:{\bigwedge}^{l-j}
\gamma^\perp\rightarrow{\bigwedge}^{l-j-1} \tau^\perp $, which are
zero maps if $\gamma$ is not a facet of $\tau$, while for
$\gamma\prec\tau$, set
$\delta_{\gamma,\tau}(w)=e_{\gamma,\tau}\lrcorner w$, where
$e_{\gamma,\tau}\in N_\QQ$ satisfies
$\tau+(-\gamma)=\RR_{\ge0}e_{\gamma,\tau}+\RR \gamma$.
 \end{defn}

\begin{rem} Compared to the definition of Ishida's complex in
\cite[Section~3.2]{o}, our definition allows the maps
$\delta_{\gamma,\tau}$ to be defined up to a rational multiple.
\end{rem}

 Ishida's complexes are also defined for star closed and star
open subsets of a fan.

\begin{defn} A subset $\Phi$  of a fan $\Sigma$ is called {\it star closed} ({\it star open}), if
$\sigma\in\Phi$ and $\sigma\prec\tau\in\Sigma$ imply $\tau\in\Phi$
(correspondingly, $\tau\in\Phi$ and $\sigma\prec\tau$ imply
$\sigma\in\Phi$).
 For such subsets,  Ishida's l-th complex of $\Phi$ is defined as
$$C^j(\Phi,l):=\bigoplus_{\begin{Sb}\gamma\in\Phi\cap\Sigma(j)\end{Sb}}
{\bigwedge}^{l-j}\gamma^\perp $$  with the coboundary homomorphism
$\delta$ as the sum of $\delta_{\gamma,\tau}$ for
$\gamma,\tau\in\Phi$.
\end{defn}

\begin{ex} Given a fan $\Sigma$ and a cone $\gamma\in\Sigma$, then
the set
$$\st_\gamma(\Sigma):=\{\tau\in\Sigma|\,\tau\succ\gamma\}$$ is
star closed in $\Sigma$. One can construct a fan
$\widetilde\st_\gamma(\Sigma)$, called the {\it star} of $\gamma$,
which consists of the images of the cones from
$\st_\gamma(\Sigma)$ in the quotient space $N_\RR/(N_\gamma)_\RR$
with the corresponding quotient lattice. Note, then, that $l$-th
Ishida's complex $C^*(\st_\gamma(\Sigma),l)$ for the subset
$\st_\gamma(\Sigma)$ is isomorphic (up to a shift) to Ishida's
complex
$C^{*-\dim\gamma}(\widetilde\st_\gamma(\Sigma),l-\dim\gamma)$. The
differential on the latter complex is induced from the first one.

\end{ex}

As in \cite[Proposition~1.8]{is}, we note the following statement.

\begin{pr} \label{p:lei} Let $\Sigma'$ be a star closed subset of a fan
$\Sigma$, equivalently $\Sigma''=\Sigma\setminus\Sigma'$ be star
open in $\Sigma$. Then there is a short exact sequence of Ishida's
complexes $$0\rightarrow C^{*}(\Sigma',l)\rightarrow
C^*(\Sigma,l)\rightarrow C^*(\Sigma'',l)\rightarrow0,$$ which
gives a long exact sequence in cohomology: $$0\rightarrow
H^0(\Sigma',l)\rightarrow H^0(\Sigma,l)\rightarrow H^0
(\Sigma'',l)\rightarrow H^1(\Sigma',l)\rightarrow
H^1(\Sigma,l)\rightarrow H^1 (\Sigma'',l)\rightarrow\cdots.$$
\end{pr}

\begin{pf} Note that if $\Sigma'$ is star closed in $\Sigma$ and $\Sigma''$ is star open
in $\Sigma$, then
 $C^{j}(\Sigma',l)$ and $C^{j}(\Sigma'',l)$ are $\QQ$-submodules  that decompose $C^j(\Sigma,l)$. Moreover,
 the differential $\delta$ for $C^{*}(\Sigma',l)$ is the same as
 $C^{*}(\Sigma,l)$,
 making it into a subcomplex.
 Similarly,
 there is a natural projection from  $C^{*}(\Sigma,l)$ onto
 $C^*(\Sigma'',l)$ compatible with the differentials.
\end{pf}


We are interested to find cohomology of the Ishida complex in the
case when the fan is a simplicial subdivision of a convex cone.
 First,  consider the simplest case, when the
fan  consists of faces of a simplicial cone.

\begin{lem}\label{l:acyc} Let $\tau$ be a simplicial cone in  $N_\RR$.
Then the higher dimensional cohomology of the complex
$$0@>>>{\bigwedge}^l
M_\QQ@>\delta>>\bigoplus_{\begin{Sb}\dim\rho=1\\
\rho\subset\tau\end{Sb}} {\bigwedge}^{l-1} \rho^\perp @>\delta>>
\cdots@>\delta>> \bigoplus_{\begin{Sb}\dim\gamma=k\\
\gamma\subset\tau\end{Sb}} {\bigwedge}^{l-k} \gamma^\perp
@>\delta>>\cdots, $$ where the differential $\delta$  is defined
as in Defintion~\ref{d:ish}, vanishes and the zeroth cohomology of
it is ${\bigwedge}^l\tau^\perp$.
\end{lem}

\begin{pf} Since the Koszul complex
$$\cdots@>\lrcorner e_\rho>>{\bigwedge}^3 M_\QQ@>\lrcorner
e_\rho>>{\bigwedge}^2 M_\QQ@> \lrcorner e_\rho>> M_\QQ @>\lrcorner
e_\rho>>\QQ@>>>0,$$ where $e_\rho$ is a $\QQ$-generator of the
1-dimensional cone $\rho$, is known to be acyclic (see
\cite[Appendix~2]{d}), we get that the zeroth cohomology of our
complex is equal to
$$\bigcap_{\rho\subset\tau}{\bigwedge}^{l}{\rho}^\perp\cong
{\bigwedge}^{l}\tau^\perp.$$

We use the induction on the dimension of the cone $\tau$ to prove
the acyclicity of the complex. The statement holds for
$\dim\tau=1$, because the map $${\bigwedge}^l M_\QQ @>\lrcorner
e_{\tau}>>{\bigwedge}^{l-1}\tau^\perp,$$ where $e_{\tau}$ is the
minimal integral generator of $\tau$, is surjective.

For $\dim\tau>1$, write $\tau=\tau'+\rho'$ where $\tau'$ and
$\rho'$ are the facet and the edge of $\tau$, respectively. Let
$\Sigma$ be the fan consisting of faces of $\tau$, and $\Sigma''$
be the fan consisting of faces of $\tau'$, then
$\Sigma'=\Sigma\setminus\Sigma''=\st_{\rho'}(\Sigma)$. So, by
Proposition~\ref{p:lei}, we get a long exact sequence in
cohomology: $$ \cdots\rightarrow
H^k(\st_{\rho'}(\Sigma),l)\rightarrow H^k(\Sigma,l)\rightarrow H^k
(\Sigma'',l)\rightarrow\cdots.$$ Note that $H^k (\Sigma'',l)$
vanishes for $k>0$ by the induction assumption. On the other hand,
$$H^k(\st_{\rho'}(\Sigma),l)\cong
H^{k-1}(\widetilde\st_{\rho'}(\Sigma),l-1).$$ Since
$\widetilde\st_{\rho'}(\Sigma)$ is a fan consisting of the faces
of the simplicial cone $\bar{\tau}$, which is the image of $\tau$
in the quotient space $N_\RR /(N_{\rho'})_\RR$, the latter
cohomology group vanishes for $k>1$, again, by the induction.
Hence, the middle term $ H^k(\Sigma,l)=0$ for $k>1$ as well. To
show that it also vanishes for $k=1$, consider the first terms of
the exact sequence:
 $$0\rightarrow
H^0(\st_{\rho'}(\Sigma),l)\rightarrow H^0(\Sigma,l)\rightarrow H^0
(\Sigma'',l)\rightarrow H^1(\st_{\rho'}(\Sigma),l)\rightarrow
H^1(\Sigma,l)\rightarrow 0  .$$ Note that
$H^0(\st_{\rho'}(\Sigma),l)=0$, since there is no zeroth term in
the corresponding complex. We already showed
$H^0(\Sigma,l)\cong{\bigwedge}^{l}\tau^\perp$ and $H^0
(\Sigma'',l)\cong{\bigwedge}^{l}{\tau'}^\perp$. Also,
$H^1(\st_{\rho'}(\Sigma),l)\cong
H^{0}(\widetilde\st_{\rho'}(\Sigma),l-1)\cong{\bigwedge}^{l-1}{\bar{\tau}}^\perp$.
To show that $H^0 (\Sigma'',l)\rightarrow
H^1(\st_{\rho'}(\Sigma),l)$ is onto, which will give the desired
result, it suffices to show that $$\dim H^0(\Sigma,l)-\dim H^0
(\Sigma'',l)+\dim H^1(\st_{\rho'}(\Sigma),l)=0.$$ But the
dimensions are the binomial coefficients $\biggl(\matrix
d-\dim\tau\\ l\endmatrix\biggr)$, $\biggl(\matrix d-\dim\tau+1\\
l\endmatrix\biggr)$ and $\biggl(\matrix d-\dim\tau\\
l-1\endmatrix\biggr)$, which satisfy the well know combinatorial
identity.
\end{pf}

The next definition is motivated by  the Chow ring of a complete
simplicial toric variety (see \cite{d}).

\begin{defn}\label{d:chow} Let $\Sigma$ be a fan in $N_\RR$ with the integral
generators of the 1-dimensional cones $e_1,\dots,e_n$. Then define
the {\it Chow} ring of $\Sigma$ as $$A
(\Sigma):={\QQ}[D_1,\dots,D_n]/(P(\Sigma)+SR(\Sigma)),$$ where
$$P(\Sigma)=\biggl\langle \sum_{i=1}^n \langle m,e_i\rangle D_i:
m\in M\biggr\rangle$$  and $$ SR(\Sigma)=\bigl\langle
D_{i_1}\cdots D_{i_k}:\{e_{i_1},\dots,e_{i_k}\}\not\subset\sigma
\text{ for all }\sigma\in\Sigma\bigr\rangle.$$ Note that $A
(\Sigma)_{*}$ is a ${\Bbb Z}$-graded algebra assuming  $\deg
D_i=1$.
\end{defn}

\begin{thm}\label{t:cohish} Let $\Sigma_\tau$ be a simplicial  subdivision of a convex
cone $\tau$ in $N_\RR$. Then   cohomology of the corresponding
$l$-th Ishida complex is $$H^k(\Sigma_\tau,l)\cong A
(\Sigma_\tau)_k\otimes {\bigwedge}^{l-k} \tau^\perp.$$
\end{thm}

\begin{pf} First, note that the $k$-th degree of the Chow ring  $A(\Sigma_\tau)$
is spanned by $D_\gamma:=\prod_{e_i\in\gamma}D_i$, for
$\gamma\in\Sigma_\tau(k)$, because of the relations
$P(\Sigma_\tau)$ and $SR(\Sigma_\tau)$. Also, the relations among
$D_\gamma$ are
\begin{equation}\label{e:relations}
\sum_{\begin{Sb}\gamma\succ\gamma'\\
\dim\gamma=\dim\gamma'+1\end{Sb}}\langle
m,e_{\gamma',\gamma}\rangle D_\gamma, \end{equation} for all
$\gamma'\in\Sigma_\tau(k-1)$, where $m\in M\cap{\gamma'}^\perp$
and $e_{\gamma',\gamma}\in N$ is a primitive generator of
$\gamma$, but not of $\gamma'$. There is a natural homomorphism
from $A(\Sigma_\tau)_k\otimes {\bigwedge}^{l-k} \tau^\perp$ to the
cohomology $H^k(\Sigma_\tau,l)$, which sends
$\sum_{\gamma\in\Sigma_\tau(k)} D_\gamma\otimes a_\gamma$, for
$a_\gamma\in{\bigwedge}^{l-k} \tau^\perp$, to $$\oplus
a_\gamma\in\ker(C^k(\Sigma_\tau,l)\rightarrow
C^{k+1}(\Sigma_\tau,l)).$$
  With this map the relations among $D_\gamma$ map to
$\im(C^{k-1}(\Sigma_\tau,l)\rightarrow C^{k}(\Sigma_\tau,l))$, so
that the homomorphism  $A(\Sigma_\tau)_k\otimes {\bigwedge}^{l-k}
\tau^\perp\rightarrow H^k(\Sigma_\tau,l)$ is well defined. We will
prove  by induction  that this homomorphism is an isomorphism in a
slightly more general situation: all  maximal cones in the fan
$\Sigma_\tau$ lie in $\tau$, have the same dimension equal to
$\dim\tau$, and the support of the fan is topologically equivalent
to a convex cone.

Suppose we are given a simplicial subdivision of $\tau$.  Pick a
cone of the dimension equal to $\dim\tau$ from this subdivision.
Then for  $\Sigma_\tau^1$ consisting of the faces of this cone,
the statement   easily follows from Lemma~\ref{l:acyc}. Let us
construct inductively, the following sequence of subfans
$\Sigma_\tau^i$ of the subdivision of $\tau$. If
$\Sigma_\tau^{i-1}$ is constructed, then $\Sigma_\tau^{i}$ is
obtained by adding to $\Sigma_\tau^{i-1}$ a new
$(\dim\tau)$-dimensional cone  from the subdivision, which is
adjacent along a facet to a cone from $\Sigma_\tau^{i-1}$, and,
also, by adding   all cones of the subdivision that have their
edges among those of $\Sigma_\tau^{i-1}$ and the new edge $\rho$
of the new cone. So, the fan $\Sigma_\tau^{i}$ has one more
1-dimensional cone, than the fan $\Sigma_\tau^{i-1}$ does. It is
clear that by doing this we get that the simplicial subdivision of
$\tau$ coincides with  $\Sigma_\tau^{i}$ for some $i$.

Now, note that  $\Sigma_\tau^{i}$ is a disjoint union of the star
open subset $\Sigma_\tau^{i-1}$ and the star closed subset
$\st_\rho(\Sigma_\tau^{i})$. So, by Proposition~\ref{p:lei} we
have a long exact sequence: $$\cdots\rightarrow
H^k(\st_\rho(\Sigma_\tau^{i}),l)\rightarrow
 H^k(\Sigma_\tau^i,l)\rightarrow H^k(\Sigma_\tau^{i-1},l)\rightarrow\cdots.$$
By induction, we can assume that $A(\Sigma_\tau^{i-1})_k\otimes
{\bigwedge}^{l-k} \tau^\perp\cong H^k(\Sigma_\tau^{i-1},l)$. We
also have  $H^k(\st_\rho(\Sigma_\tau^{i}),l)\cong
H^{k-1}(\widetilde\st_\rho(\Sigma_\tau^{i}),l-1)\cong
A(\widetilde\st_\rho(\Sigma_\tau^{i}) )_{k-1}\otimes
{\bigwedge}^{l-k} \bar\tau^\perp $, where the induction is applied
to the fan $\widetilde\st_\rho(\Sigma_\tau^{i})$ lying in the
image $\bar\tau$ of the cone $\tau$ in the quotient space
$(N/N_\rho)_\RR$. Using the description of $A(\Sigma_\tau^{i})$ in
terms of $D_\gamma$, one can easily show that there is a natural
exact sequence: $$0\rightarrow
A(\widetilde\st_\rho(\Sigma_\tau^{i}) )_{k-1}\rightarrow
A(\Sigma_\tau^{i})_k\rightarrow
A(\Sigma_\tau^{i-1})_k\rightarrow0.$$ Hence, we get a commutative
diagram  with exact rows and  with the right- and left-hand
vertical maps being isomorphisms: $$\minCDarrowwidth0.5cm
\begin{CD}
0@>>> A(\widetilde\st_\rho(\Sigma_\tau^{i}) )_{k-1}\otimes V
 @>>> A(\Sigma_\tau^{i})_k \otimes V
@>>> A(\Sigma_\tau^{i-1})_k\otimes V @>>>0\\ @. @VVV @VVV @VVV
@.\\0@>>> H^k(\st_\rho(\Sigma_\tau^{i}),l)@>>>
H^k(\Sigma_\tau^i,l)@>>> H^k(\Sigma_\tau^{i-1},l)@>>>0,
\end{CD}$$
where $V={\bigwedge}^{l-k} \bar\tau^\perp=
{\bigwedge}^{l-k}\tau^\perp$. By the  5-lemma (see \cite{e}), the
middle vertical map is also isomorphism.
\end{pf}

\begin{cor} \label{c:dim} Let $\Sigma_\tau$ be a simplicial  subdivision of a convex
cone $\tau$ in $N_\RR$. Then  $$\dim H^k(\Sigma_\tau,l)=
\biggl(\matrix d-\dim\tau\\ l-k\endmatrix\biggr)\cdot\sum_{j=0}^k
\biggl(\matrix \dim\tau-j\\ k-j\endmatrix\biggr)
(-1)^{k-j}\cdot\#\Sigma_\tau(j).$$
\end{cor}

\begin{pf} First, we guess the formula based on the answer for
complete simplicial fans (see \cite[Corollary~4.2]{o2}). Then
check that it also holds for a fan $\Sigma_\tau$ consisting of
faces of a simplicial cone $\tau$ and  prove it by the induction
using the exact sequence from the proof of the above theorem.
\end{pf}

\begin{rem} One can consider Ishida's complexes of $\ZZ$-modules. However, for a
fan with singular cones, the corresponding cohomology of Ishida's
complexes may have torsion. \end{rem}

\section{Cohomology of rational forms and a vanishing
theorem.}\label{s:two}

Our goal in this section is to compute the cohomology
$H^{k}(\ps,\Omega_{\ps}^{l}(X))$ for all $k$, $l$  and a semiample
divisor $X$ on a complete simplicial toric variety $\ps$. While
not all of these spaces vanish for $k>0$ as it was in the case of
an ample divisor, we discover that some of them do vanish. The
result for a trivial divisor  gives the  well-known description of
the cohomology of a complete simplicial toric variety.

 Before we can state our results, let us review some futher notation.
 Let ${\bf P}_{\Sigma}$ be a $d$-dimensional complete  toric variety associated
with  the fan  $\Sigma$ in $N_\RR$.  We denote by ${\bf T}_\sigma$
 a torus corresponding to the cone $\sigma\in\Sigma$ and by $V(\sigma)$
 the closure of ${\bf T}_\sigma$ in $\ps$.
Also, $D_1,\dots,D_n$ are the torus invariant irreducible divisors
in ${\bf P}_\Sigma$, corresponding to the  primitive generators
 $e_1,\dots,e_n$  of the 1-dimensional cones.
The polynomial ring  $S=S(\ps)={\Bbb C}[x_1,\dots,x_n]$ is called
the homogeneous coordinate ring of the toric variety $\ps$. A
torus invariant Weil divisor $D=\sum_{i=1}^{n}a_iD_i$ on the
complete toric variety gives rise to a convex polytope
$$\Delta_D=\{m\in M_{\Bbb R}:\langle m,e_i\rangle\geq-a_i \text{
for all } i\}\subset M_{\Bbb R}.$$ There is also a support
function $\psi_D : N_{\Bbb R}\rightarrow\Bbb R$ which is linear on
each cone $\sigma\in\Sigma$ and $\psi_D(e_i)=\langle m_\sigma,
e_i\rangle=-a_i$ for all $e_i\in\sigma$ and  some $m_\sigma\in M$.

We will work with semiample divisors which  are conveniently
classified by the following definition.

\begin{defn} \cite{m} A semiample Cartier divisor $D$ (i.e.,
${\cal O}_{\ps}(D)$ is generated by global sections) on a complete
toric variety $\ps$ is called {\it $i$-semiample} if  the
Kodaira-Iitaka dimension $\kappa(D):=\dim\phi_D(\ps)=i$, where
$\phi_D:\ps@>>>{\Bbb P}(H^0(\ps,{\cal O}_{\ps}(D)))$ is the
rational map defined by the sections of the line bundle ${\cal
O}_{\ps}(D)$.
\end{defn}

These divisors satisfy  the  property:

\begin{thm}\label{t:fun} \cite{m}
Let $[D]\in A_{d-1}(\ps)$ be an $i$-semiample divisor class on a
complete toric variety $\bold P_\Sigma$ of dimension $d$. Then,
there exists a unique complete toric variety $\psd$ with a
surjective morphism $\pi:{\bf P}_\Sigma@>>>{\bf P}_{\Sigma_D}$,
arising {from} a surjective homomorphism of lattices
$\tilde\pi:N@>>>N_D$ which maps the fan $\Sigma$ into $\Sigma_D$,
such that $\pi^*[Y]=[D]$ for some ample divisor $Y$ on $\psd$.
Moreover,  $\dim\psd=i$, and, for a torus invariant $D$, the fan
$\Sigma_D$ in $(N_D)_{\Bbb R}:=N_{\Bbb R}/N'_{\Bbb R}$, where
$N'=\{v\in N:\psi_D(-v)=-\psi_D(v)\}$ is a sublattice of $N$ and
$\psi_D$ is the support function of $D$,
 is the normal fan of $\Delta_D$, which lies in $(M_D)_{\Bbb R}$,
where $M_D:={N'}^\perp\cap M$.
\end{thm}

To perform the calculations of the cohomology of the twisted
sheaves we need to use  Ishida's complex of sheaves which is a
resolution of $\Omega_{\ps}^{l}$ (see \cite[Section~3.2]{o}). Let
us recall its construction. Define $${\cal
I}^{k,l}_\ps:=\bigoplus_{\gamma\in\Sigma(k)}\Omega_{V(\gamma)}^{l-k}(\log
D(\gamma)), $$ for $0\le k\le l$, where $D(\gamma):=\sum_{\gamma'}
V({\gamma'})$ is the anticanonical divisor on $V(\gamma)$,
 and set ${\cal I}^{k,l}_\ps:=0$ for $k<0$ or $k>l$.
Then, there are natural morphisms $${\cal
I}^{k,l}_\ps@>\delta>>{\cal I}^{k+1,l}_\ps,$$ induced by the
Poincar\'e residue maps $$
R_{{\gamma'},\gamma}:\Omega_{V(\gamma)}^{l-k}(\log D(\gamma))
@>>>\Omega_{V(\gamma')}^{l-k-1} (\log D({\gamma'}))$$
 for the component $V({\gamma'})$ of the divisor $D(\gamma)$ on $V(\gamma)$,
if $\gamma\in\Sigma(k)$ is a face of ${\gamma'}\in\Sigma(k+1)$.
 If
$\gamma\in\Sigma(k)$ is not a face of ${\gamma'}\in\Sigma(k+1)$,
then the corresponding map $R_{{\gamma'},\gamma}=0$.

By \cite[Theorem~3.6]{o}, the following sequence is exact
$$0@>>>\Omega_{\ps}^{l}@>>>{\cal I}^{0,l}_\ps@>\delta>>{\cal
I}^{1,l}_\ps @>\delta>>\cdots@>\delta>>{\cal I}^{l,l}_\ps@>>>0,$$
since $\ps$ is simplicial. Twisting this sequence by ${\cal
O}_\ps(X)$  gives another exact sequence
\begin{equation}\label{e:rest}
0@>>>\Omega_{\ps}^{l}(X)@>>>{\cal I}^{0,l}_\ps(X)@>\delta>>{\cal
I}^{1,l}_\ps(X) @>\delta>>\cdots@>\delta>>{\cal
I}^{l,l}_\ps(X)@>>>0,
\end{equation}
which is a resolution of the twisted sheaf. Here, ${\cal
I}^{k,l}_\ps(X):={\cal I}^{k,l}_\ps\otimes{\cal O}_\ps(X)$. We
will now use this resolution to compute the cohomology of the
twisted sheaves.

\begin{thm} \label{t:cohrf} Let $\beta=[X]\in A_{d-1}(\ps)$ be a semiample divisor class
on a complete simplicial toric variety $\ps$, and let
$\pi:\ps@>>>\psx$ be the associated canonical contraction. Then
$$H^{k}(\ps,\Omega_{\ps}^{l}(X))\cong\bigoplus\begin{Sb}\sigma\in\Sigma_X
\end{Sb}
\bigl(S/\langle x_j:\,\tilde\pi(\rho_j)\subset\sigma
\rangle\bigr)_{\beta-\beta_0+\beta_1^\sigma}\otimes
A^\sigma(\Sigma)_k\otimes
{\bigwedge}^{l-k}(M_X\cap\sigma^\perp),$$ where
$\beta_0=\deg(\prod_{j=1}^n x_j)$,
$\beta_1^\sigma=\deg(\prod_{\tilde\pi(\rho_j)\subset\sigma}x_j)$,
and $$A^\sigma(\Sigma):={\Bbb
Q}[D_1,\dots,D_n]/(P(\Sigma)+SR(\Sigma)+ \langle
D_j:\tilde\pi(\rho_j)\not\subset\sigma\rangle)$$ is defined
similar to Definition~\ref{d:chow}.
\end{thm}

\begin{pf}
Since
\begin{equation}\label{e:iden}
\Omega_{V(\gamma)}^{l-k}(\log D(\gamma))\simeq {\cal
O}_{V(\gamma)}\otimes{\bigwedge}^{l-k}(M\cap\gamma^\perp)
\end{equation}
(see \cite[Corollary~3.2]{o}), the sheaf ${\cal I}^{k,l}_\ps(X)$
is a direct sum of the semiample
 sheaves ${\cal O}_{V(\gamma)}(X)$.
By the vanishing of the higher dimensional cohomology of a
semiample sheaf (see \cite[Corollary~7.3]{d}), it follows
 that the resolution (\ref{e:rest}) is acyclic, whence the cohomology
of the twisted sheaf $\Omega_{\ps}^{l}(X)$ can be computed as the
cohomology of the complex of global sections (see
\cite[Proposition~4.3]{ii}): $$0@>>>H^0(\ps,{\cal
I}^{0,l}_\ps(X))@>\delta>> H^0(\ps,{\cal I}^{1,l}_\ps(X))
@>\delta>>\cdots@>\delta>>H^0(\ps,{\cal I}^{l,l}_\ps(X))@>>>0.$$
By the identifications (\ref{e:iden}) and the isomorphisms (see
\cite{c1}) $$H^0(V(\gamma),{\cal O}_{V(\gamma)}(X))\cong
S(V(\gamma))_{\beta^\gamma},$$ where $\sigma\in\Sigma_X$ is the
smallest cone containing the image of $\gamma$, this complex can
be rewritten as
\begin{multline}\label{e:comp}
0@>>>S_\beta\otimes{\bigwedge}^l
M@>\delta>>\bigoplus_{\rho\in\Sigma(1)}
S(V(\rho))_{\beta^\rho}\otimes{\bigwedge}^{l-1} (M\cap\rho^\perp)
@>\delta>>
\\
\cdots@>\delta>>\bigoplus_{\gamma\in\Sigma(k)}
S(V(\gamma))_{\beta^\gamma}\otimes{\bigwedge}^{l-k}
(M\cap\gamma^\perp)
@>\delta>>\cdots@>\delta>>\bigoplus_{\gamma\in\Sigma(l)}S(V(\gamma))_{\beta^\gamma}@>>>0,
\end{multline}
where the map $\delta$ sends $g\otimes w\in
S(V(\gamma))_{\beta^\gamma}\otimes{\bigwedge}^{l-k}
(M\cap\gamma^\perp)$, for $\gamma\in\Sigma(k)$, to the direct sum
of $g_{\gamma'}\otimes (w\lrcorner e_{{\gamma'},\gamma})\in
S(V({\gamma'}))_{\beta^{\gamma'}}\otimes{\bigwedge}^{l-k-1}
(M\cap{\gamma'}^\perp)$, for
$\gamma\subset{\gamma'}\in\Sigma(k+1)$, such that $g_{\gamma'}$ is
the image of $g$ induced by the restriction $H^0(V(\gamma),{\cal
O}_{V(\gamma)}(X))@>>>H^0(V({\gamma'}),{\cal
O}_{V({\gamma'})}(X))$, and $e_{{\gamma'},\gamma}$ is the minimal
integral generator of the cone $\rho\in\Sigma(1)$ contained in
${\gamma'}$ but not in $\gamma$.

Next, notice that the monomials in $S_\beta$ naturally  correspond
to the lattice points of the polytope $\Delta_X$, by \cite{c1}.
With respect to this identification, the monomials in
$$S(V(\gamma))_{\beta^\gamma}\cong \bigl(S/\langle
x_j:\,\tilde\pi(\rho_j)\subset\sigma\rangle\bigr)_\beta$$ (see
\cite{m}) are the lattice points in the face of $\Delta_X$
corresponding to the minimal cone $\sigma\in\Sigma_X$ containing
the image of $\gamma$. The natural grading of
$S(V(\gamma))_{\beta^\gamma}$, for $\gamma\in\Sigma$, by the
lattice points of $\Delta_X$ induces a grading on the sequence
(\ref{e:comp}), and it is not difficult to see that the maps
$\delta$ in  (\ref{e:comp}) respect this grading. A monomial in
$S_\beta$ can be uniquely written as
$\prod_{\tilde\pi(\rho_j)\not\subset\sigma}x_j$ times a monomial
in $\bigl(S/\langle x_j:\,\tilde\pi(\rho_j)\subset\sigma
\rangle\bigr)_{\beta-\beta_0+\beta_1^\sigma}$, where
$\sigma\in\Sigma_X$ corresponds to the minimal face of $\Delta_X$
containing the lattice point associated to the monomial. {From}
here, it follows  that the $k$-th cohomology of (\ref{e:comp}) is
isomorphic to  the direct sum, by $\sigma\in\Sigma_X$, of the
tensor products of the complex spaces $\bigl(S/\langle
x_j:\,\tilde\pi(\rho_j)\subset\sigma
\rangle\bigr)_{\beta-\beta_0+\beta_1^\sigma}$ and the $k$-th
cohomology of the complex
\begin{equation*}
0@>>>{\bigwedge}^l
M@>\delta>>\bigoplus_{\begin{Sb}\rho\in\Sigma(1)\\
\tilde\pi(\rho)\subset\sigma\end{Sb}} {\bigwedge}^{l-1}
(M\cap\rho^\perp) @>\delta>> \cdots@>\delta>>
\bigoplus_{\begin{Sb}\gamma\in\Sigma(k)\\
\tilde\pi(\gamma)\subset\sigma\end{Sb}} {\bigwedge}^{l-k}
(M\cap\gamma^\perp) @>\delta>>\cdots.
\end{equation*}
Recognize that this is the Ishida complex of $\ZZ$-modules for a
simplicial subdivision of the convex cone $\tilde\pi^{-1}(\sigma)$
induced by the fan $\Sigma$. Since we tensor the cohomology groups
of this complex with complex spaces, we can discard torsion and
use Theorem~\ref{t:cohish}. After noting
$M\cap(\tilde\pi^{-1}(\sigma))^\perp=M_X\cap\sigma^\perp$ the
result easily follows.
\end{pf}

As a consequence of the above theorem, we get the following
vanishing result, which is an analogue of  the
Bott-Steenbrink-Danilov vanishing theorem (see \cite{bflm}).

\begin{thm}\label{t:vanish}  Let $X$ be an $i$-semiample divisor
on a complete simplicial toric variety $\ps$. Then
$$H^{k}(\ps,\Omega_{\ps}^{l}(X))=0$$ if $k>l$ or $l>k+i$.
\end{thm}

\begin{pf}
If $k>l$, then ${\bigwedge}^{l-k}(M_X\cap\sigma^\perp)=0$ in
Theorem~\ref{t:cohrf}. One can also obtain
$H^{k}(\ps,\Omega_{\ps}^{l}(X))=0$ in this case directly from the
complex of global sections of sheaves arising from the exact
sequence (\ref{e:rest}), since that complex does not have nonzero
terms after $l$.

 If $l>k+i$ then $l-k>i$, while the rank of $M_X\cap\sigma^\perp$
 is no more than $i$. So, ${\bigwedge}^{l-k}(M_X\cap\sigma^\perp)$
 vanishes again for all $\sigma\in\Sigma_X$.
\end{pf}
\begin{rem}
When $X$ is a trivial divisor (i.e., $i=0$), this theorem gives
the vanishing part of the cohomology of the complete simplicial
toric variety $\ps$: $H^k(\ps,\Omega_\ps^l)=0$ for $k\ne l$.
Moreover, by Theorem~\ref{t:cohrf}, we get
$H^k(\ps,\Omega_\ps^k)\cong \CC\otimes A(\Sigma)_k$, which is the
Danilov-Jurkevich description of the cohomology of a complete
simplicial toric variety with coefficients in $\CC$ (see
\cite{d}). We should also remark that Theorem~\ref{t:vanish} can
not be extended to  arbitrary complete toric varieties, because in
the case when $X$ is trivial we   get
$H^{k}(\ps,\Omega_{\ps}^{l})$, which may be nontrivial for $k<l$
(see \cite{ad}). However, for $k>l$ and a complete toric variety
$\ps$, cohomology $H^{k}(\ps,\Omega_{\ps}^{l})=0$ by \cite{d}.
 It will be
interesting to see if the vanishing result of
Theorem~\ref{t:vanish} holds   for complete toric varieties when
$k>l$.
\end{rem}

We can also deduce  a generalization of the Kodaira vanishing
theorem for toric varieties, which was proved by a different
method in \cite{mu} and \cite{BBo} (see also \cite{cd}).

\begin{thm} Let $X$ be an $i$-semiample divisor
on a $d$-dimensional complete simplicial toric variety $\ps$. Then
$H^{k}(\ps,\Omega_{\ps}^{d}(X))=0$  for $k\ne d-i$.
\end{thm}

\begin{pf} If $l=d$ in Theorem~\ref{t:cohrf}, then
${\bigwedge}^{l-k}(M_X\cap\sigma^\perp)\ne 0$ for $d-k\le
i-\dim\sigma$. For such $\sigma$, we have
$\dim\tilde\pi^{-1}(\sigma)=\dim\sigma+d-i\le k$. By the relations
in $A^\sigma(\Sigma)_k$, it is clear that this $\QQ$-module
vanishes for $\dim\tilde\pi^{-1}(\sigma)< k$. Let us show that in
the case of $\dim\tilde\pi^{-1}(\sigma)=k $, we also have
$A^\sigma(\Sigma)_k=0$. Indeed, the maximal cones $\gamma$ of
dimension $k$ that subdivide the convex cone
$\dim\tilde\pi^{-1}(\sigma)$ correspond to the generators
$D_\gamma$ of $A^\sigma(\Sigma)_k$. By the relations
(\ref{e:relations}) that come from facets of the maximal cones,
all $D_\gamma$ are multiples of each other. Moreover, they are all
zero, if $\sigma\ne\{0\}$, since  we also have the relation coming
from the facet of a maximal cone, which lies on the boundary of
$\tilde\pi^{-1}(\sigma)$. But $\sigma\ne\{0\}$ because
$\dim\tilde\pi^{-1}(\sigma)=k\ne d-i $ while
$\dim\tilde\pi^{-1}(\{0\})=d-i$.
\end{pf}

\begin{cor}\label{c:formula}  Let   $ X $ be an $i$-semiample torus invariant divisor
on a $d$-dimensional complete simplicial toric variety $\ps$ as in
Theorem~\ref{t:cohrf}. Then $$\dim H^{k}(\ps,\Omega_{\ps}^{l}(X))=
\sum_{\Gamma\prec\Delta_X} l^*(\Gamma)\biggl(\matrix  \dim\Gamma\\
l-k\endmatrix\biggr)\cdot\sum_{j=0}^k \biggl(\matrix
d-\dim\Gamma-j\\ k-j\endmatrix\biggr)
(-1)^{k-j}\cdot\#\Sigma_{\sigma_\Gamma}(j),$$ where the sum is by
faces of the polytope $\Delta_X$,   ${\sigma_\Gamma}\in\Sigma_X$
corresponds to $\Gamma$,
$\Sigma_{\sigma_\Gamma}=\{\tau\in\Sigma|\,\tilde\pi(\tau)\subset\sigma_\Gamma\}$
is a subfan of $\Sigma$, and $l^*(\Gamma)$ denotes the number of
interior lattice points inside $\Gamma$.
\end{cor}

\begin{pf} Cones $\sigma$ of $\Sigma_X$ correspond to the faces $\Gamma$ of the polytope $\Delta_X$
and monomials in $\bigl(S/\langle
x_j:\,\tilde\pi(\rho_j)\subset\sigma
\rangle\bigr)_{\beta-\beta_0+\beta_1^\sigma}$ correspond to
monomials in $S_\beta$ that are divisible by
$\prod_{\tilde\pi(\rho_j)\not\subset \sigma}x_j$ and not divisible
by $x_j$ for $\tilde\pi(\rho_j) \subset \sigma$. Such monomials
are in one-to-one correspondence with the interior lattice points
of $\Gamma$. Next, note that $A^\sigma(\Sigma) \simeq
A(\Sigma_{\tilde\pi^{-1}(\sigma)})$, where
$\Sigma_{\tilde\pi^{-1}(\sigma)}$ is the simplicial subdivision of
$\tilde\pi^{-1}(\sigma)$ induced by $\Sigma$. Since
 $ M_X\cap\sigma^\perp = M\cap\tilde\pi^{-1}(\sigma)^\perp $,
 Corollary~\ref{c:dim} and Theorem~\ref{t:cohrf} give us the
 formula.
\end{pf}

\begin{ex} E. Materov in \cite{mat} calculated the Bott formula
for ample divisors on complete simplicial toric varieties. Let
  $X$   be ample on $\ps$, then $X$ is $d$-semiample and
  $\Sigma_X=\Sigma$ consists of simplicial cones. For $k=0$ in Corollary \ref{c:formula},  we get
$$\dim H^{0}(\ps,\Omega_{\ps}^{l}(X))=
\sum_{\Gamma\prec\Delta_X} l^*(\Gamma)\biggl(\matrix  \dim\Gamma\\
l\endmatrix\biggr).$$ For $k>0$, the combinatorial identity
$$\sum_{j=0}^k\biggl(\matrix
m-j\\ k-j\endmatrix\biggr) (-1)^{k-j}\biggl(\matrix m \\
 j\endmatrix\biggr)=0,$$ which follows from
 $$\sum_{j=0}^k \biggl(\matrix m \\
 j\endmatrix\biggr) t^j(1-t)^{m-j}=(t+1-t)^m=1,$$
 implies
 $\dim H^{k}(\ps,\Omega_{\ps}^{l}(X))=0$.

If $X$ is trivial on $\ps$, then $\Sigma_X$ and $\Delta_X$ are
just  points. We get
$$\dim H^{k}(\ps,\Omega_{\ps}^{l})=\Biggl\{\matrix
  \sum_{j=0}^k \biggl(\matrix
d -j\\ k-j\endmatrix\biggr) (-1)^{k-j}\cdot\#\Sigma (j), & k=l,
\\
 0, \hspace{4cm} & k\ne l.\endmatrix$$
 These are the formulas in Theorems 2.14 and  3.6 in \cite{mat}.
\end{ex}

 The next
thing we will do is to describe $ H^{k}(\ps,\Omega_{\ps}^{l}(X))$
via \v Cech cohomology. The toric variety $\ps$ has an affine open
cover ${\cal U}=\{U_\tau\}$, where
$$U_\tau=\{x\in\ps|\,\prod_{e_i\notin\tau}x_i\ne0\}$$ and
$\tau\in\Sigma(d)$ have the maximal dimension. This   cover
induces an affine open cover on all subvarieties
$V(\gamma)\subset\ps$ as well. Using the notation in \cite{bc},
let $m\in M$ correspond to the differential form
$\sum_{i=1}^n\langle m,e_i\rangle\frac{\dd x_i}{x_i}$.

\begin{pr}\label{p:cocyle} Let $ X $ be a semiample divisor
on a complete simplicial toric variety $\ps$ defined by $f=0$ for
$f\in S_\beta$. Then, under the isomorphism of
Theorem~\ref{t:cohrf} and the natural isomorphism $\check
H^{k}({\cal U}, \Omega_{\ps}^{l}(X))\cong
H^{k}(\ps,\Omega_{\ps}^{l}(X))$, we have that $A\otimes
D_{i_1}\cdots D_{i_k}\otimes \omega$, where  $A\in S
_{\beta-\beta_0+\beta_1^\sigma}$, $\gamma=\langle
e_{i_1},\dots,e_{i_k}\rangle\in \Sigma$,
$\tilde\pi(\gamma)\subset\sigma\in\Sigma_X$ and $\omega
\in{\bigwedge}^{l-k}(M_X\cap\sigma^\perp)$,
 correspond to the \v Cech
cocycle $$\biggl\{\frac{A}{f}\prod_{\tilde\pi(e_i)\notin\sigma}x_i
(m_{\tau_1}^{ i_1}-m_{\tau_0}^{ i_1})
\wedge(m_{\tau_2}^{i_{2}}-m_{\tau_1}^{ i_{2}})\wedge\cdots\wedge
(m_{\tau_k}^{i_k}-m_{\tau_{k-1}}^{
i_k})\wedge\omega\biggr\}_{\tau_0\dots\tau_k}$$ with
$m_{\tau}^{i}=0$ if $e_i\notin\tau$, and, for $e_i\in\tau$,
$m_{\tau}^{ i}$ is determined by $\langle m_{\tau}^{
i},e_i\rangle=1$ and $\langle m_{\tau}^{ i},e_j\rangle=0$ for
$e_j\in\tau$, $j\ne i$.

\end{pr}

\begin{pf} We have the resolution (\ref{e:rest}) of the sheaf of
rational forms  $\Omega_{\ps}^{l}(X)$. When dealing with an
acyclic resolution of sheaves, one can apply the following
standard trick. Introduce auxiliary sheaves ${\cal K}_j$ as
kernels of the morphisms ${\cal I}^{j,l}_\ps(X)@>\delta>>{\cal
I}^{j+1,l}_\ps(X)$. Then, by the exactness of (\ref{e:rest}), we
get short exact sequences of sheaves $$0\rightarrow{\cal
K}_j\rightarrow{\cal I}^{j,l}_\ps(X)\rightarrow {\cal K}_{j+1},$$
which give rise to the long exact sequences in cohomology:
\begin{multline*}
\cdots\rightarrow H^{k-j}(\ps,{\cal I}^{j,l}_\ps(X) )\rightarrow
H^{k-j}(\ps,{\cal K}_j)\rightarrow
\\
H^{k-j+1}(\ps,{\cal K}_{j-1})\rightarrow H^{k-j+1}(\ps,{\cal
I}^{j,l}_\ps(X) ) \rightarrow\cdots.
\end{multline*}
 Since higher dimensional cohomology
of ${\cal I}^{j,l}_\ps(X)$ vanishes, the connecting homomorphisms
are isomorphisms for $k-j>0$ and epimorphisms for $k-j=0$. Thus,
we get isomorphisms:\begin{align*}H^0(\ps,{\cal K}_k)/\im
H^0(\ps,{\cal I}^{k-1,l}_\ps(X))&\cong H^1(\ps,{\cal
K}_{k-1})\cong H^2(\ps,{\cal K}_{k-2})\cong\cdots
\\
&\cong H^{k-1}(\ps,{\cal K}_1)\cong
H^{k}(\ps,\Omega_{\ps}^{l}(X)).\end{align*}

Now, ${A}(\prod_{\tilde\pi(e_i)\notin\sigma}x_i)\omega/f$ is a
global section of $\Omega_{V(\gamma)}^{l-k}(\log D(\gamma))$,
which is a subsheaf of ${\cal I}^{k,l}_\ps(X)$. Moreover, one can
easily check that its restriction by the homomorphism
$$H^0(\ps,{\cal I}^{k,l}_\ps(X))@>\delta>>H^0(\ps,{\cal
I}^{k+1,l}_\ps(X))$$ vanishes, whence
${A}(\prod_{\tilde\pi(e_i)\notin\sigma}x_i)\omega/f\in
H^0(\ps,{\cal K}_k)$. In \v Cech cohomology this section is
represented by the cocycle
\begin{equation}\label{e:cocycle1}
\biggl\{
{A}\biggl(\prod_{\tilde\pi(e_i)\notin\sigma}x_i\biggr)\omega/f\biggr\}_\tau\in
\check H^0({\cal U}|_{V(\gamma)},\Omega_{V(\gamma)}^{l-k}(\log
D(\gamma))).
\end{equation} To find its
image by the connecting homomorphism $\check H^0(\ps,{\cal
K}_k)\rightarrow \check H^0(\ps,{\cal K}_{k-1})$, we can use the
commutative diagram: $$ \minCDarrowwidth0.4cm
\begin{CD}
0 @>>> C^{1}({\cal U} ,{\cal K}_{k-1} ) @>>> C^{1}({\cal U} ,{\cal
I}^{k-1,l}_\ps(X) ) @>>> C^{1}({\cal U} ,{\cal K}_{k})@>>>0
\\ @. @AAA @AAA @AAA   \\ 0 @>>>C^{0}({\cal U} , {\cal K}_{k-1})
@>>> C^{0}({\cal U} ,{\cal I}^{k-1,l}_\ps(X)) @>>> C^{0}({\cal U}
, {\cal K}_{k})@>>>0.
\end{CD}
$$ The homomorphism $C^{0}({\cal U} ,{\cal I}^{k-1,l}_\ps(X)) @>>>
C^{0}({\cal U} , {\cal K}_{k})$ is induced by the Poincar\'e
residue maps $$R_{{\gamma},\gamma'}:C^0({\cal
U}|_{V(\gamma')},\Omega_{V(\gamma')}^{l-k+1}(\log D(\gamma'))
\rightarrow C^0({\cal U}|_{V(\gamma)},\Omega_{V(\gamma)}^{l-k}
(\log D({\gamma})),$$ where $\gamma'$ is a facet of
$\gamma\in\Sigma(k)$. The cocycle (\ref{e:cocycle1}) has a lift by
this homomorphism to the cochain
$$\biggl\{\frac{A}{f}\prod_{\tilde\pi(e_i)\notin\sigma}x_i
m_{\tau}^{ i_k} \wedge\omega\biggr\}_{\tau }\in C^0({\cal
U}|_{V(\gamma')},\Omega_{V(\gamma')}^{l-k+1}(\log D(\gamma')),$$
where $\gamma'=\langle e_{i_1},\dots,e_{i_{k-1}}\rangle\in
\Sigma(k-1)$, which can also be thought as the cochain in
$C^{0}({\cal U} ,{\cal I}^{k-1,l}_\ps(X))$. Applying the \v Cech
coboundary to this cochain, we get the cocycle
 $$\biggl\{\frac{A}{f}\prod_{\tilde\pi(e_i)\notin\sigma}x_i
(m_{\tau_1}^{i_k}-m_{\tau_{0}}^{
i_k})\wedge\omega\biggr\}_{\tau_0\tau_1}.$$

Repeating the above procedure and using the commutative diagrams
$$ \minCDarrowwidth0.4cm
\begin{CD}
0 @>>> C^{k-j+1}({\cal U} ,{\cal K}_{j-1} ) @>>> C^{k-j+1}({\cal
U} ,{\cal I}^{j-1,l}_\ps(X) ) @>>> C^{k-j+1}({\cal U} ,{\cal
K}_{k})@>>>0
\\ @. @AAA @AAA @AAA   \\ 0 @>>>C^{k-j}({\cal U} , {\cal K}_{j-1})
@>>> C^{k-j}({\cal U} ,{\cal I}^{j-1,l}_\ps(X)) @>>> C^{k-j}({\cal
U} , {\cal K}_{j})@>>>0,
\end{CD}
$$ one can check that the image of the cocycle (\ref{e:cocycle1})
by the sequence of homomorphisms $\check H^{k-j}(\ps,{\cal
K}_j)\rightarrow \check H^{k-j+1}(\ps,{\cal K}_{j-1})$ is
precisely  the cocycle
$$\biggl\{\frac{A}{f}\prod_{\tilde\pi(e_i)\notin\sigma}x_i
(m_{\tau_1}^{ i_1}-m_{\tau_0}^{ i_1})
\wedge(m_{\tau_2}^{i_{2}}-m_{\tau_1}^{ i_{2}})\wedge\cdots\wedge
(m_{\tau_k}^{i_k}-m_{\tau_{k-1}}^{
i_k})\wedge\omega\biggr\}_{\tau_0\dots\tau_k}$$ from $\check
H^{k}({\cal U}, \Omega_{\ps}^{l}(X))$.
\end{pf}

\begin{ex} Let us  apply Proposition~\ref{p:cocyle} in the case
 when $\ps$ is the projective space $\PP^n$,
and $X$ is a trivial divisor. If $x_1,\dots,x_{n+1}$ are the
homogeneous coordinates on $\PP^n$, then the open cover
 ${\cal U}$ is given by the open sets
 $U_\tau=\{x\in\PP^n|\,x_j\ne0\}$, where $\tau=\langle
 e_1,\dots,\widehat{e_j},\dots,e_{n+1}\rangle$, and
 $m_\tau^i=\frac{\dd x_i}{x_i}-\frac{\dd x_j}{x_j}$.
Hence, the following cocycles represent
$H^k(\PP^n,\Omega^k_{\PP^n})\cong \CC$: $$\biggl\{ \biggl(
\frac{\dd x_{j_1}}{x_{j_1}}-\frac{\dd x_{j_0}}{x_{j_0}}\biggr)
\wedge\biggl(\frac{\dd x_{j_2}}{x_{j_2}}-\frac{\dd
x_{j_1}}{x_{j_1}}\biggr)\wedge\cdots\wedge \biggl(\frac{\dd
x_{j_k}}{x_{j_k}}-\frac{\dd
x_{j_{k-1}}}{x_{j_{k-1}}}\biggr)\biggr\}_{j_0\dots j_k}.$$

\end{ex}

\end{document}